\documentclass[a4paper,12pt]{article}
\usepackage[utf8]{inputenc}
\usepackage[english]{babel}
\usepackage{enumitem}
\usepackage{color}
\usepackage{amsmath}
\usepackage{amssymb}
\usepackage{amsthm}
\usepackage{amsfonts}
\usepackage{epsfig}
\usepackage{subcaption}
\usepackage{caption}
\usepackage{graphicx}
\usepackage{epstopdf}
\usepackage{multirow}
\usepackage{upgreek}
\usepackage{placeins}
\usepackage{mathrsfs}
\usepackage{float}%
\usepackage{booktabs}%??????\toprule??\midrule??\bottomrule
\usepackage{multirow}
\usepackage{pifont}

\usepackage{caption}
\newtheorem{dfn}{Definition}[section]
\newtheorem{lem}{Lemma}[section]

\newtheorem{alg}{Algorithm}[section]
\newtheorem{theorem}{Theorem}[section]
\theoremstyle{definition}
\newtheorem{asm}{Assumption}[section]

\newtheorem{rem}{Remark}[section]

\newtheorem{fact}{Fact}[section]
\usepackage{tabularx}  % 自适应宽度
\usepackage{booktabs}

\DeclareMathOperator{\diag}{diag}
\DeclareMathOperator*{\zer}{zer}
\DeclareMathOperator*{\Fix}{Fix}
\DeclareMathOperator*{\gra}{gra}

\DeclareMathOperator*{\Id}{Id}
\DeclareMathOperator*{\range}{range}

\DeclareMathOperator*{\dom}{dom}

\DeclareMathOperator*{\slt}{slt}
\allowdisplaybreaks
\begin{document}
	\title{\textbf{Frugal forward-backward splitting methods with deviations}}
	\author{ \sc \small Yongyu Fu$^a${\thanks{ email: fyy000611@163.com}},\,\,
		Haowen Zheng$^b${\thanks{email: zhw20251122@163.com}},\,\,
		Qiao-Li Dong$^a${\thanks{Corresponding author. email: dongql@lsec.cc.ac.cn}},\,\,
Xiaolong Qin$^c${\thanks{Corresponding author. email: qxlxajh@163.com}}\,\,\,and\,
		Jing Zhao$^a${\thanks{email: zhaojing200103@163.com}}\\
		\small $^a$College of Science,  Civil Aviation University of China, Tianjin 300300, China,\\
		\small $^b$School of Light Industry and Engineering, South China University of Technology,\\ 	\small  Guangzhou 510640, China,\\
\small $^c$Department of Mathematics, Hangzhou Normal University, Hangzhou 311121, China\\
	}
	\date{}
	\maketitle
	
	\begin{abstract}
		
		{
			The deviation vectors provide additional degrees of freedom and effectively enhance the flexibility of algorithms.
			In the literature, the iterative schemes with deviations are constructed and their convergence analyses
			are performed on an inefficient, algorithm-by-algorithm basis.
			In this paper, we  address these by providing a general framework of frugal forward-backward splitting methods with deviations for finding zeros in the sum of a finite number of maximally monotone operators and  cocoercive operators. Our framework encompasses the Douglas--Rachford splitting method with deviations.
			A unified weak convergence analysis is made  under mild conditions.  Numerical experiments on  Markowitz portfolio optimization problem are given to demonstrate the effectiveness of deviations.	
		}
	\end{abstract}
	

	\section{ Introduction }
	In this paper, we focus on the structured monotone inclusion problem in a real Hilbert space $\mathcal{H}$ which is to find $x\in \mathcal{H}$  such that
	\begin{equation}
		\label{ABC}
		0  \in \left ( \sum_{i=1}^{n}F_i+\sum_{i=1}^{m}B_i   \right )(x),
	\end{equation}
	where $F_1,\dots,F_n : \mathcal{H} \rightrightarrows \mathcal{H}$ are   maximally monotone and   $B_i: \mathcal{H}\to \mathcal{H} $ are $\frac{1}{L_i}$-cocoercive,  $i=1,\dots,m$.  The abstraction given by  (\ref{ABC}) covers many fundamental problems arising in mathematical optimization involving objectives with finite-sum structure, such as  composite minimization, structured saddle-point problems and structured variational inequalities \cite[Examples 1, 2, 3]{ring}.
	
	Based on the splitting structural form of the problem (\ref{ABC}), researchers have proposed some  forward-backward-type splitting  methods \cite{Moursi,Chen}. These methods use only vector addition, scalar multiplication, the resolvents of the  maximally monotone operators and direct evaluations of the cocoercive  operators.
	Frugal forward-backward splitting methods with the minimal lifting in the sense of Ryu \cite{Ryu} have recently received  much attention \cite{ring,fran,Frugal,Bred}.  They  calculate the resolvent of each maximally monotone operator and each cocoercive operator only once in each iteration.  Moreover, these methods minimize storage requirements.
	In the recent seminal  work, Morin et al. \cite{Frugal}  showed that the minimal lifting number $n-1$ of the frugal forward-backward splitting methods for \eqref{ABC} is achievable only if  the first and last operator evaluations are resolvent evaluations.
	
	For solving the problem \eqref{ABC} with $n\ge2$ and $m=n-1$, some scholars designed several frugal forward-backward type splitting methods with the minimal lifting.
	Arag\'{o}n-Artacho et al. \cite{ring} firstly extended the   techniques from  \cite{Malit} to  propose a distributed forward-backward splitting method.
	Based on the product space trick \cite{Pierra} and  the degenerate preconditioned
	proximal point algorithms, Bredies et al. \cite{Bred}
	developed  parallel and sequential forward-backward splitting methods by generalizing Davis--Yin splitting method.
	
	To systematically investigate the frugal forward-backward splitting methods with  minimal lifting and efficiently analyze their convergence,  there has been increasing interest  recently in constructing algorithmic frameworks for solving \eqref{ABC}. As far as we know, three frameworks have been constructed.
	Arag\'{o}n-Artacho et al. \cite{fran} firstly extended  methodologies in \cite{Bre} and introduced a graph-based framework by incorporating an additional spanning subgraph to model the cocoercive operators' role for  the problem \eqref{ABC} with $n\ge2$ and $m=n-1$.
	Then \AA kerman et al. \cite{Anton}  proposed a more general  framework  for  the problem \eqref{ABC} with $n\ge1$ and $m\ge1$ as follows.
	\begin{alg}\label{AlgG}
		\hrule\vskip 0.5mm
		\rm\noindent\rm%\small{{A frugal forward-backward  algorithmic framework with deviations.}}
		\vskip 0.5mm
		\hrule
		
		\vskip 1mm
		\noindent
		%\noindent\textbf{Pick:} The parameters $\eta, \theta$,  $\{\gamma_k\}_{k\in \mathbb{N}} $ and   $\{\xi_k\}_{k\in \mathbb{N}}$ satisfying Assumption \ref{assumption2}, and  the matrices $M\in \mathbb{R}^{n\times (n-1)},  S\in \mathbb{R}^{n\times n}, C\in \mathbb{R}^{n\times m} $ and $Q\in \mathbb{R}^{m\times n}$ satisfying                                                                                                                                                                                                                                                                                                                                                                                                                                                                                                                                                                                                                                                                                                                                                                                                      Assumption \ref{asm1}. \\ \textbf{Let:} $(d_1,\dots,d_n)=2(\diag(S))^{\odot(-1)}$.\\
		1: \textbf{Input:} $\gamma\in(0,1)$, $\mathbf{z}^0=(z_1^0,\dots,z_{n-1}^0)\in \mathcal{H}^{n-1}$.\\ %$\mathbf{u}^0=(0,\dots,0)\in \mathcal{H}^{m}$ and $\mathbf{v}^0=(0,\dots,0)\in \mathcal{H}^{n-1}$. \\
		2: \textbf{for}  $k=0,1,2,...$  \textbf{do} \\
		3:\ \ \ \ \  \textbf{for} $i=1,\dots,n$ \textbf{do}\\
		$$
		x_i^{k}=J_{d_i F_i}\left({-d_i\sum_{j=1}^{i-1}}S_{ij}x_j^{k}+d_i{ \sum_{j=1}^{n-1}}M_{ij}(z_j^k)-d_i\sum_{j=1}^{m}C_{ij} B_{j}\left(  { \sum_{h=1}^{i-1}} Q_{jh}x_h^k \right)\right).
		$$
		4:\ \ \ \ \  \textbf{end for}\\
		5:  $z_i^{k+1}=z_i^k-\gamma{\textstyle \sum_{i=1}^{n}} M_{ij}x_i^k, \quad\forall i \in[1,n-1].$\\
		6: \textbf{end for}
		\vskip 1mm
		
		\hrule
		
		\hspace*{\fill}
	\end{alg}
	\noindent In Algorithm \ref{AlgG},  $S,M,C,Q$ are matrices satisfying some conditions, and $(d_1,\dots,d_n)$ $=2(\diag(S))^{\odot(-1)}$.
	It was  verified in \cite{Anton} that Algorithm \ref{AlgG} encompasses Davis--Yin splitting method and  the graph-based algorithmic framework in \cite{fran}. At almost the same time, Dao et al. \cite{minh} adopted approaches similar to those in \cite{Anton} and  devised a general framework for solving \eqref{ABC} with $n\ge2$ and $ m= n-1$ by introducing  matrices satisfying specific assumptions. It is worth mentioning that the stepsize in \cite{fran,minh}  depends on $\min_{i}\{L_i\}$ while that in \cite{Anton} involves each $L_i$, $i=1,\ldots,m$. Therefore, the latter is generally larger than the former.

	Recently, the variants of the methods through deviations have attracted much attention since  Sadeghi et al.  \cite{Sedeghi} first proposed a forward-backward splitting algorithm with deviations by incorporating the deviation vectors into each iteration. To increase the flexibility in the  process of learning to optimize \cite{L2O} while guaranteeing convergence, a deviation-based proximal gradient method was developed in \cite{Ban}. It learns the entire update function and then uses it in an  optimization solver.  Qin et al. \cite{Qin} also investigated a  Douglas--Rachford (DR) splitting algorithm with deviations whose scheme is given as
	\begin{equation}
		\label{Qin}
		\left\{
		\aligned
		&x_1^k=J_{\gamma F_1}(z^k+v^k),\\
		&x_2^{k}=J_{\gamma F_2}(2x_1^k-(z^k+v^k)),\\
		&z^{k+1}=z^k-\lambda_k(x_1^k-x_2^k),
		\endaligned
		\right.
	\end{equation}
	where $\lambda_k\in(0,2)$, $\gamma>0$ and $\{v^k\}$ is  deviation vectors.  In order to ensure the convergence of the algorithm  (\ref{Qin}), the deviation vectors $\{v^k\}$ are required to satisfy the norm inequality
	\begin{equation}\label{norm}
		\frac{\lambda_{k+1}}{2-\lambda_{k+1}}  \| v^{k+1}\|^2 \le\xi_k\lambda_k(2-\lambda_k) \| x_2^k-x_1^k+\frac{1}{2-\lambda_k}v^k  \|^2 ,
	\end{equation}
	where $\xi_k\in[0,1)$. 	
	In addition,  there have emerged  some other methods with deviations, such as Davis--Yin splitting method with deviations \cite{Hu-Davis-Yin} and forward-backward-half forward splitting algorithm with deviations as well \cite{Qin-FBHF}. Despite many similarities in the convergence proofs of these algorithms (e.g.,
	the norm inequalities of deviations are introduced to obtain a nonincreasing sequence), their convergence analyses are  performed on an inefficient, algorithm-by-algorithm basis in the literature.
	
	In this work, we first propose  a frugal forward-backward splitting algorithmic framework with deviations by incorporating deviation vectors into Algorithm \ref{AlgG}.  Our framework  avoids the drawback of adding deviations to algorithms one by one. Then we derive a norm inequality  to ensure the convergence. With the aid of the norm inequality, we  establish the  weak convergence results of the proposed algorithmic framework.

	The structure of the paper is as follows. In Section 2, we introduce the notations and main concepts. In Section 3, we  propose a frugal forward-backward splitting algorithmic framework with deviations  and give the weak convergence analysis. In Section 4, we present numerical experiments.
	
	\section{Preliminaries}
	Let $\mathbb{N}=\{0,1,\dots\}$ be the set of natural numbers and $\mathbb{N}_+=\{1,2,\dots\}$ be the set of non-zero natural numbers. Throughout this paper,  we denote by $\mathcal{H}$ a real Hilbert space with inner product  $\langle \cdot, \cdot\rangle$ and induced norm $\|\cdot\|=\sqrt{\langle \cdot, \cdot\rangle}$.  Denote by $\Id$  the identity operator on $\mathcal{H}$. %A bounded, self-adjoint  operator $H: \mathcal{H}\to \mathcal{H}$ is said to be strongly positive if there exists some $c>0$ such that  $\langle x,Hx \rangle\ge c \| x  \|^2 $ for all $x\in\mathcal{H}$.
	Given a linear, self-adjoint, and strongly  positive  operator $H: \mathcal{H}\to \mathcal{H}$, we define $\langle x,y  \rangle_H =\langle x,Hy  \rangle $  and  $ \| x  \|_H=\sqrt{ \langle x,Hx  \rangle }  $ for $\forall x,y \in\mathcal{H}$.     We use $\omega_w(x^k) = \{x : \exists x^{k_j}\rightharpoonup x\}$ to denote the weak $\omega$-limit set of the sequence $\{x^k\}$.
	
	\begin{fact}
		For all $a, b, c, d\in\mathcal{H}$, there holds	
		\begin{equation}
			\label{ceq}
			2\left \langle a-b,c-d \right \rangle =\left \| a-d \right \| ^2+\left \| b-c \right \|^2-\left \| a-c \right \|^2-\left \| b-d \right \|  ^2.
		\end{equation}
	\end{fact}	
	
	Let   $T: \mathcal{H} \rightarrow {\mathcal{H}}$ be an operator. Denote by $\Fix T$ be the set of fixed points of $T,$ i.e., $\Fix T=\{x\in \mathcal{H} \,:\, x=T(x)\}.$
	\begin{dfn}
		\rm
		An operator  $T: \mathcal{H} \rightarrow {\mathcal{H}}$ is said to be
		\begin{itemize}
			\item[(i)]  $\beta$-Lipschitz continuous, if there exists a constant $ \beta> 0$, such that
			$$
			\|T(x)-T(y)\|\leq  \beta\|x-y\|,\quad\forall x, y\in \mathcal{H},
			$$
			and nonexpansive if $\beta=1$.
			\item[(ii)] $\frac{1}{\beta}$-cocoercive, if there exists a constant $ \beta> 0$, such that
			$$
			\langle T(x)-T(y),x-y \rangle \geq  \frac{1}{\beta}\|T(x)-T(y)\|^2,\quad\forall x, y\in \mathcal{H}.
			$$
		\end{itemize}
	\end{dfn}
	\noindent
	By  Cauchy--Schwarz inequality, a $\frac{1}{\beta}$-cocoercive operator is $\beta$-Lipschitz continuous.
	
	Given a set-valued operator $F:\mathcal{H}\rightrightarrows\mathcal{H}$, the domain, the range,  the graph and the zeros of $F$ are respectively denoted by
	$\dom F =\{x\in\mathcal{H} : F(x)\neq\varnothing\}$,  $\range F =\{u\in\mathcal{H}: u\in F(x)\ \hbox{for}\ \forall x\in\dom F\}$,
	$\gra F =\{(x,u)\in\mathcal{H}\times\mathcal{H} : u\in F(x)\}$  and
	$\zer F =\{x\in\mathcal{H} : 0\in F(x)\}$.
	The inverse operator of $F$, denoted by $F^{-1}$, is defined through $x\in F^{-1}(u)\Leftrightarrow u\in F(x)$.

	\begin{dfn}{\rm(\cite[Definition 20.1 and Definition 20.20]{BC2011})}
		{\rm
			A set-valued operator $F: \mathcal{H} \rightrightarrows\mathcal{H}$ is said to be
			\begin{itemize}
				\item[(i)]  monotone if $\langle u-v,x-y \rangle \geq 0$, $\forall (x,u),(y,v) \in \gra F $.
				\item[(ii)]  maximally monotone if there exists no monotone operator $B: \mathcal{H} \rightrightarrows\mathcal{H}$ such that $\gra B$ properly contains $\gra F,$ i.e., for every $(x,u) \in \mathcal{H}\times\mathcal{H}$
				$$
				(x,u) \in \gra F  \ \ \Leftrightarrow \ \ \langle u-v , x-y \rangle \geq 0,  \ \  \forall(y,v)\in \gra F.
				$$
			\end{itemize}
		}
	\end{dfn}
	
	Given an operator $F : \mathcal{H}\rightrightarrows\mathcal{H}$, the resolvent of $F$ with parameter $\lambda>0$ is  denoted by $J_{\lambda F} =(\Id+\lambda F)^{-1}$. From \cite[Corollary 23.11]{BC2011}, if the operator $F$ is  maximally  monotone, then $J_{\lambda F}$ is single-valued and 1-cocoercive.

	For a matrix $P\in\mathbb{R}^{n\times m}$, we denote by $P_{ij}$ its $(i,j)$ component. The transpose of matrix $P$ is  denoted as  $P^{\top}$.  When $m=n,$ we denote by $\slt(P)\in\mathbb{R}^{n\times n}$ the strictly  lower triangular matrix extracted from $P$, and  by   $\diag(P)\in\mathbb{R}^n$ the vector   extracted  from the main diagonal of $P$.  With a mild overload of notation, for $w\in\mathbb{R}^n$, we denote by $\diag(w)\in\mathbb{R}^{n\times n}$ the diagonal matrix with the diagonal being $w$. Given a vector $g = (g_1,\dots, g_n)$ with all non-zero elements, the Hadamard inverse of $g$, denoted by $g^{\odot (-1)}$, is the vector defined element-wise by $(g^{\odot (-1)})_i = \frac{1}{g_i} $ for  $i=1,\dots,n$.
	
	Given a matrix $P\in \mathbb{R} ^{n\times m}$,  we denote the Kronecker  product of $P$ and $\Id$ by
	$$
	\mathbf{P} =P\otimes \Id=\begin{bmatrix}
		P_{11}\Id&  P_{12}\Id & \cdots   & P_{1m}\Id\\
		P_{21}\Id & P_{22}\Id  & \cdots &  P_{2m}\Id\\
		\vdots  & \vdots  &\ddots   &\vdots  \\
		P_{n1}\Id &  P_{n2}\Id &\cdots   & P_{nm}\Id
	\end{bmatrix}.
	$$
	Note that $\mathbf{P}$ is a bounded linear operator from $\mathcal{H}^m$ to $\mathcal{H}^n$.
	
	\begin{dfn}\rm($(m,n)$-nondecreasing vector).
		Let $n\in \mathbb{N}_+$, $m\in \mathbb{N}$. We say that  $$A=(A_1,A_2,\dots,A_n)\in\{0,1,\dots,m\}^n$$ is an $(m,n)$-nondecreasing vector if $A_1=0,$ $A_n=m$ and $A_i\le A_{i+1}$ for each $i\in[1,n-1]$.
	\end{dfn}
	Let $n\in \mathbb{N}_+$, $m\in \mathbb{N}$ and $A$ be an $(m,n)$-nondecreasing  vector. Denote by $\mathcal{S}(A)$  the set of matrices $R\in \mathbb{R}^{n\times m}$ that have a  staircase structure w.r.t. $A$, i.e.,
	$$R_{ij}
	=0\ \ \hbox{for}\ \hbox{all} \ i=1,\dots,n\ \hbox{and}\ j>A_i.$$
	Conversely, $\mathcal{S}^c(A)$ denotes the set of matrices $R\in\mathbb{R}^{n\times m}$ that have a complement staircase structure w.r.t. $A$, i.e., $R_{ij}=0$ for all $i=1,\dots,n$ and $j\le A_i$.
	
	Next we use the concepts of staircase and complement staircase to define what is called a causal pair of matrices.
	\begin{dfn}\label{causal}\rm
		(Causal pair of matrices). A pair of matrices $C,Q^{\top}\in \mathbb{R}^{n\times m}$ is said to be causal if there exists an $(m,n)$-nondecreasing  vector $A$ such that
		$$C\in\mathcal{S}(A)\ \hbox{and}\ Q^{\top}\in\mathcal{S}^c(A).$$
	\end{dfn}
	\begin{lem}\label{lle}{\rm(\cite[Lemma 2.47]{BC2011})}
		Let $\Theta$ be a nonempty set of $\mathcal{H}$, and $\{x^k\}$ be a sequence in $\mathcal{H}$.  Assume that the following conditions hold:\vskip 1mm
		\begin{itemize}
			\item[{\rm(i)}]  for every $x \in \Theta$, $\lim_{k \rightarrow \infty}\|x^{k}-x\|$ exists;
			\noindent
			\item[{\rm(ii)}] every weak sequential cluster point of $\{x^k\}$ belongs to $\Theta$, i.e., $\omega_w(x^k)\subseteq \Theta$.
		\end{itemize}
		Then the sequence $\{x^k\}$ converges weakly to a point in $\Theta$.
	\end{lem}
	
	%	\begin{rem}\rm The concept of causal pairs of matrices captures the full class of matrices $C$ and $Q^{\top}\in\mathbb{R}^{n\times m}$ such that, given $\mathbf{B}:=\diag(B_1,\dots,B_m)$ with $B_i:\mathcal{H}\to\mathcal{H}$, the term $C\mathbf{B}(Q\mathbf{x}^{k})$ in Algorithm \ref{Alg1} can be computed only using a single evaluation of each operator $B_i$ and simple algebraic operations. And $C\mathbf{B}Q$ is strictly lower triangular for all choices of $B_j$, which is consistent with the ordering defined by $A$.
		%	\end{rem}
	
	\section{Main results}
	
	In this section, we first present a frugal forward-backward splitting algorithmic framework with deviations. Then the weak convergence of the proposed framework is established.
	
	Throughout this section, the solution set
	of the problem \eqref{ABC} is assumed to be nonempty, i.e., $\zer( {\textstyle \sum_{i=1}^{n}}  F_i+{\textstyle \sum_{i=1}^{m}} B_i)\neq \emptyset $. Define an operator   $\mathbf{F} : \mathcal{H}^n\rightrightarrows \mathcal{H}^n$ by $\mathbf{F(x)} =\huge (  F_1(x_1),\dots,$ $F_n(x_n)\huge )$ for all $\mathbf{x}=(x_1,\dots,x_n)\in\mathcal{H}^n$. It follows from \cite[Proposition 20.23]{BC2011} that $\mathbf{F}$ is  maximally monotone  and its resolvent $J_{\mathbf{F}}: \mathcal{H}^n\to \mathcal{H}^n$ is given by $J_{\mathbf{F}}=(J_{F_1}, \dots,J_{F_n})$. Similarly define an operator $\mathbf{B}:\mathcal{H}^{m}\to \mathcal{H}^{m}$ by $\mathbf{B(y)} =(B_1(y_1),\dots,B_{m}(y_m))$ for all $\mathbf{y}=(y_1,\dots,y_m)\in\mathcal{H}^m$. Let  $L=\diag(\frac{1}{L_1},\dots,\frac{1}{L_m})$, where $\frac{1}{L_i}$ is the cocoercive constant of $B_i$ for $i=1,\dots,m$. From \cite[Lemma 7.2]{Frugal} , it follows that
	\begin{equation}\label{BB}
		\langle \mathbf{B}(\mathbf{x})-\mathbf{B}(\mathbf{y}),\mathbf{x}-\mathbf{y}\rangle\ge\|\mathbf{B}
		(\mathbf{x})-\mathbf{B}(\mathbf{y})\|_{\mathbf{L}},
		\quad\forall \mathbf{x},\mathbf{y}\in \mathcal{H}^{m}.
	\end{equation}
	
	\subsection{Algorithm}
	
	Before giving our algorithm,  assume that the parameters $\theta$,  $\{\gamma_k\} $ and  $\{\xi_k\}$ satisfy the following conditions.
	\begin{asm}
		\label{assumption2}
		\rm
		Let $\epsilon$ be a sufficiently small positive constant. For all $k\in\mathbb{N}$, assume that  the follows hold:
		\begin{itemize}
			\item[{\rm(i)}] $\theta \in (0,+\infty)$;
			\item[{\rm(ii)}]    $\epsilon\le \gamma_k \le 1-\epsilon$;
			\item[{\rm(iii)}]   $0\le \xi _k\le1-\epsilon$.
		\end{itemize}
	\end{asm}	
	\noindent	
	We also need some matrices such that the following assumptions hold.
	\begin{asm}
		\label{asm11}
		The  matrices $M\in \mathbb{R}^{n\times (n-1)},  S\in \mathbb{R}^{n\times n}, C\in \mathbb{R}^{n\times m} $ and $Q\in \mathbb{R}^{m\times n}$ satisfy the following properties:
		\begin{itemize}
			\item[{\rm(a)}]$\ker M^{\top}\subseteq \mathbb{R}e_n$, where $e_n=(1,\dots,1)^{\top}\in \mathbb{R}^n$.
			\item[{\rm(b)}]$C$, $Q^{\top}$  are causal  and   $C^{\top}e_n=Qe_n=e_m,$ where $e_m=(1,\dots,1)^{\top}\in\mathbb{R}^m.$
			\item[{\rm(c)}]$S$ is symmetric    such that $e_n^{\top}Se_n=0$ and $S-MM^{\top}-\frac{1}{2}(1+\frac{1}{\theta})(C^{\top}-Q)^{\top}{L^{-1}}(C^{\top}-Q)\succeq 0$.
		\end{itemize}
	\end{asm}
	\noindent

	Now we  present a frugal forward-backward splitting algorithmic framework with deviations for solving  \eqref{ABC} with $n\ge1$ and $m\ge0$.
	\begin{alg}\label{Alg12}
		\hrule\vskip 0.5mm
		\noindent\rm\small{{A frugal forward-backward  algorithmic framework with deviations.}}
		\vskip 0.5mm
		\hrule
		
		\vskip 1mm
		\noindent
		\noindent\textbf{Pick:} The parameters $\theta$,  $\{\gamma_k\}$ and   $\{\xi_k\}$ satisfying Assumption \ref{assumption2}, and  the matrices $M\in \mathbb{R}^{n\times (n-1)},  S\in \mathbb{R}^{n\times n}, C\in \mathbb{R}^{n\times m} $ and $Q\in \mathbb{R}^{m\times n}$ satisfying                                                                                                                                                                                                                                                                                                                                                                                                                                                                                                                                                                                                                                                                                                                                                                                                      Assumption \ref{asm11}. \\ \textbf{Let:} $(d_1,\dots,d_n)=2(\diag(S))^{\odot(-1)}$.\\
		1: \textbf{Input:} $\mathbf{z}^0=(z_1^0,\dots,z_{n-1}^0)\in \mathcal{H}^{n-1}$, $\mathbf{u}^0=(0,\dots,0)\in \mathcal{H}^{m}$ and $\mathbf{v}^0=(0,\dots,0)\in \mathcal{H}^{n-1}$. \\
		2: \textbf{for}  $k=0,1,2,...$  \textbf{do} \\
		3:\ \ \ \ \  \textbf{for} $i=1,\dots,n$ \textbf{do}\\
		$$x_i^{k}=J_{d_i F_i}\left({-d_i \sum_{j=1}^{i-1}}S_{ij}x_j^{k}+d_i{ \sum_{j=1}^{n-1}}M_{ij}(z_j^k+v_j^k)-d_i\sum_{j=1}^{m}C_{ij} B_{j}\left(  { \sum_{h=1}^{i-1}} Q_{jh}x_h^k +u_{j}^k\right)\right).$$
		4: \ \ \ \ \  \textbf{end for}\\
		5:   $z_i^{k+1}=z_i^k-\gamma_k{\textstyle \sum_{i=1}^{n}} M_{ij}x_i^k, \quad\forall i \in[1,n-1].$\\
		6:  Choose $\mathbf{u}^{k+1}$ and $\mathbf{v}^{k+1}$ such that the following inequality is satisfied:\\
		\begin{equation}\label{tiaojian}
			\frac{\gamma_{k+1}}{1-\gamma_{k+1}} \| \mathbf{v}^{k+1}  \|^2+\frac{\gamma_{k+1} (1+\theta )}{2} \| \mathbf{u}^{k+1}  \|_{\mathbf{L}^{-1}} ^2\le \xi _kl^2_k ,
		\end{equation}
		\ \ \ \	where
		\begin{equation}\label{ln}
			l_k^2=\frac{1-\gamma_k}{\gamma_k}\| \mathbf{z}^{k+1}-\mathbf{z}^k +\frac{\gamma_k}{1-\gamma_k}\mathbf{v}^k \|^2.
		\end{equation}
		7: \textbf{end for}
		\vskip 1mm
		
		\hrule
		
		\hspace*{\fill}
	\end{alg}
	Now we discuss the relation of Algorithm \ref{Alg12} and the Douglas--Rachford splitting method with deviations \eqref{Qin}-\eqref{norm}.
	Let $n=2$, $m=0$,
	\begin{equation*}
		\label{DRd}
		\aligned
		M=\sqrt{\frac{2}{\gamma}}\begin{bmatrix}1
			\\
			-1
		\end{bmatrix},\  S=\begin{bmatrix}
			\frac{2}{\gamma}& -\frac{2}{\gamma} \\
			-\frac{2}{\gamma}&  \frac{2}{\gamma}
		\end{bmatrix},
		\endaligned
	\end{equation*}
	and $d_1=d_2=\gamma$. Then Algorithm \ref{Alg12} becomes a Douglas--Rachford splitting method with deviations:
	\begin{equation}
		\label{DRd1}
		\left\{
		\aligned
		&x_1^k=J_{\gamma F_1}(\sqrt{2\gamma}z^k+\sqrt{2\gamma}v^k),\\
		&x_2^{k}=J_{\gamma F_2}(2x_1^k-(\sqrt{2\gamma}z^k+\sqrt{2\gamma}v^k)),\\
		&z^{k+1}=z^k-\gamma_k\sqrt{\frac{2}{\gamma}}(x_1^k-x_2^k).
		\endaligned
		\right.
	\end{equation}
	To compare with that given by \eqref{Qin}-\eqref{norm},
	set $\bar{z}^k=\sqrt{2\gamma}z^k$, $\bar{v}^k=\sqrt{2\gamma}v^k$
	and $\lambda_k=2\gamma_k$. Then \eqref{DRd1} can be reformulated as
	\begin{equation*}
		\left\{
		\aligned
		&x_1^k=J_{\gamma F_1}(\bar{z}^k+\bar{v}^k),\\
		&x_2^{k}=J_{\gamma F_2}(2x_1^k-(\bar{z}^k+\bar{v}^k)),\\
		&\bar{z}^{k+1}=\bar{z}^k-\lambda_k(x_1^k-x_2^k),
		\endaligned
		\right.
	\end{equation*}
	which  obviously equals  to \eqref{Qin}.
	From \eqref{tiaojian} and \eqref{ln}, it follows that  the deviation vector $\bar v^k$  satisfies the norm inequality
	\begin{equation}\label{pic}
		\frac{\lambda_{k+1}}{2-\lambda_{k+1}}\| \bar v^{k+1}\|^2 \le\xi_k\lambda_k(2-\lambda_k) \| x_2^{k}-x_1^k +\frac{1}{2-\lambda_k}\bar{v}^k \|^2.
	\end{equation}
	It is clear that the deviation inequalities \eqref{pic} and \eqref{norm} are equivalent. Therefore, our framework includes the Douglas--Rachford splitting method with deviations \eqref{Qin}-\eqref{norm} as a special case.

	For the convenience of the convergence analysis, we  introduce an operator $T:\mathcal{H}^{n-1}\to \mathcal{H}^{n-1}$  of the form
	\begin{equation}\label{T}
		T(\mathbf{z})= \mathbf{z}-\gamma\mathbf{M}^*\mathbf{x},\ \ \mathbf{x}=J_{ \mathbf{ D F}}\mathbf{(DNx+DMz-D CB(Qx))},
	\end{equation}
	where $N=-\slt(S)$ and $D=\diag(d_1,\dots,d_n)$.
	Note that	Algorithm \ref{Alg12} can be recast as
	\begin{equation}\label{framework2}
		\mathbf{z}^{k+1}=\mathbf{z}^k-\gamma_k \mathbf{M}^*\mathbf{x}^k,\
		\mathbf{x}^k=J_{ \mathbf{DF}}(\mathbf{DNx}^k+\mathbf{DM}(\mathbf{z}^k+\mathbf{v}^k)-\mathbf{DCB}
		(\mathbf{Q}\mathbf{x}^k+\mathbf{u}^k)).
	\end{equation}

	\begin{rem}
		\label{rem11}
		\begin{itemize}
			\item[{\rm(1)}]
			A simple choice of $S$ satisfying Assumption \ref{asm11}(c) is $S=MM^{\top}+PP^{\top}+\frac{1}{2}(1+\frac{1}{\theta})(C^{\top}-Q)^{\top}{L^{-1}}(C^{\top}-Q)$, where $P\in \mathbb{R}^{n\times (n-1)}$ is such that  $\mathbb{R}e_n\subseteq \ker P^{\top}$.  
			In Algorithm \ref{Alg12}, $d_i$, $i=1,\ldots,n$ is  the stepsize. Obviously, to obtain large stepsizes, $P=0$ is  the optimal option for $S$.
			
			\item[{\rm(2)}]
			By $N=-\slt(S)$ and $D=\diag(d_1,\dots,d_n)$, we get $S=2D^{-1}-N-N^{\top}$.
			According to Assumption \ref{asm11}(c),	the matrices $M\in \mathbb{R}^{n\times (n-1)},$ $N\in \mathbb{R}^{n\times n}$,  $D\in \mathbb{R}^{n\times n} $,  $C\in \mathbb{R}^{n\times m}$ and $Q\in \mathbb{R}^{m\times n}$ satisfy the following properties:
			\begin{itemize}
				\item[{\rm(a)}] $e_n^{\top}(D^{-1}-N)e_n=0$;
				\item[{\rm(b)}]$2D^{-1}-N-N^{\top}\succeq MM^{\top}+\frac{1}{2}(1+\frac{1}{\theta})(C^{\top}-Q)^{\top}{L^{{-1}}}(C^{\top}-Q)$.
			\end{itemize}
			
			\item[{\rm(3)}]
			The concept of causal pairs of matrices captures the full class of matrices $C$, $Q^{\top}\in\mathbb{R}^{n\times m}$ such that the term $\mathbf{CB(Q}\mathbf{x}^{k}+\mathbf{u}^{k})$ in  \eqref{framework2} can be computed only using a single evaluation of each operator $B_i$ and simple algebraic operations. Furthermore, $\mathbf{CB(Q})$ is strictly lower triangular for all choices of $B_i$ according to Definition \ref{causal}.
			
		\end{itemize}
	\end{rem}

	The next lemma gives the relation of the fixed point set of $T$ and the  zeros of ${\textstyle \sum_{i=1}^{n}}  F_i+{\textstyle \sum_{i=1}^{m}} B_i$.
	One can easily show it by employing arguments which are similar to those used in the proof of  \cite[Lemma 3.5]{minh} and \cite[Lemma 3.1]{Matt}.
	\begin{lem}
		\label{lem11}
		(Fixed points and zeros) Suppose that Assumption \ref{asm11}(a) and (b) hold, and let
		$$
		\aligned
		\Omega=\{(\mathbf{z},x)\in \mathcal{H}^{n-1}\times \mathcal{H}:\ &\mathbf{x}=J_{ \mathbf{DF}}(\mathbf{DNx+DMz- DCB(Qx)}) \,\,\hbox{where}\\
		&\mathbf{x}=(x,\dots,x)\in \mathcal{H}^{n} \}.
		\endaligned
		$$ Then the following assertions hold.
		\begin{itemize}
			\item[{\rm(a)}] If $\mathbf{z}\in \Fix T$, then there exists $x\in \mathcal{H}$ such that $(\mathbf{z},x)\in \Omega$.
			\item[{\rm(b)}] If $x\in \zer \left (  {\textstyle \sum_{i=1}^{n}} F_i+ {\textstyle \sum_{i=1}^{m}B_i}  \right ) $, then there exists $\mathbf{z}\in \mathcal{H}^{n-1}$ such that $(\mathbf{z},x)\in \Omega$.
			\item[{\rm(c)}] If $(\mathbf{z},x)\in \Omega$, then $\mathbf{z}\in \Fix T $ and $x\in \zer \left (  {\textstyle \sum_{i=1}^{n}} F_i+ {\textstyle \sum_{i=1}^{m}B_i}  \right ) $.
			Consequently, $$\Fix T \ne \emptyset \Leftrightarrow \Omega \ne \emptyset \Leftrightarrow \zer \left (  {\textstyle \sum_{i=1}^{n}} F_i+ {\textstyle \sum_{i=1}^{m}B_i}  \right )\ne \emptyset.$$
		\end{itemize}
	\end{lem}
	
	\subsection{Convergence analysis}
	Now we present a lemma which is key for the  convergence analysis of Algorithm \ref{Alg12}.
	
	\begin{lem}\label{lem22}
		Let $F_1,\dots,F_n:\mathcal{H}\rightrightarrows \mathcal{H}$ be maximally monotone, and   let $B_i: \mathcal{H}\to \mathcal{H} $ be $\frac{1}{L_i}$-cocoercive,  $i=1,\dots,m$  with  $\zer({\textstyle \sum_{i=1}^{n}} F_i+{\textstyle \sum_{i=1}^{m}} B_i)\neq \emptyset $. Suppose that Assumption \ref{assumption2} and  Assumption \ref{asm11}(c) hold. Let $\{\mathbf{z}^k\}\subseteq\mathcal{H}^{n-1}$ and  $\{\mathbf{x}^k\}\subseteq\mathcal{H}^{n}$  be sequences generated by  Algorithm \ref{Alg12}, and  $\{\mathbf{u}^k\}\subseteq\mathcal{H}^{m}$ and $\{\mathbf{v}^k\}\subseteq\mathcal{H}^{n-1}$  satisfy \eqref{tiaojian}. Let $\mathbf{z}=(z_1,\dots,z_{n-1})\in\mathcal{H}^{n-1}$ be an arbitrary point in  $\Fix T$. Then, for all $k\in\mathbb{N}$,  we have
		\begin{equation}\label{lem1}
			\| \mathbf{z}^{k+1}-\mathbf{z}  \|^2+l_k^2
			\le  \| \mathbf{z}^k-\mathbf{z}  \|^2+ \frac{\gamma_k}{1-\gamma_k}\|\mathbf{v}^k\|^2+\frac{\gamma_k (1+\theta )}{2}  \| \mathbf{u}^k  \| _{\mathbf{L}^{{-1}}}^2.
		\end{equation}
		Furthermore,
		\begin{equation}
			\label{dfd}
			\| \mathbf{z}^{k+1}-\mathbf{z}  \|^2+l_k^2\le  \| \mathbf{z}^{k}-\mathbf{z}  \|^2+\xi _{k-1}l_{k-1}^2
		\end{equation}
		holds for all $k\in\mathbb{N}_+$.
	\end{lem}
	\begin{proof}
		Let  $\mathbf{x}=J_{\mathbf{ DF}}(\mathbf{y})$ with $\mathbf{y}=\mathbf{DMz+DNx- DCB(Qx)}$ and $\mathbf{x}^k=J_{\mathbf{DF}}(\mathbf{y}^k)$ with  $\mathbf{y}^k=\mathbf{DM}(\mathbf{z}^k+\mathbf{v}^k)+\mathbf{DNx}^k-\mathbf{ DCB(Q}\mathbf{x}^k+\mathbf{u}^k)$. Since $\mathbf{D}^{-1}\mathbf{(y-x)}\in\mathbf{F}(\mathbf{x})$ and $\mathbf{D}^{-1}\mathbf{(y}^k-\mathbf{x}^k)\in\mathbf{F}(\mathbf{x}^k)$, the monotonicity of $\mathbf{F}$ gives
		\begin{equation}\label{chang}
			\aligned
			0  \leq&\langle\mathbf{x}-\mathbf{x}^k,\mathbf{D}^{-1}(\mathbf{y}-\mathbf{x})-\mathbf{D}^{-1}(\mathbf{y}^k-\mathbf{x}^k)\rangle \\
			=&\langle\mathbf{x}-\mathbf{x}^k,\mathbf{D}^{-1}(\mathbf{DMz}+\mathbf{DNx}- \mathbf{DCB(Qx)}-\mathbf{x})\rangle \\
			&-\langle\mathbf{x}-\mathbf{x}^k,\mathbf{D}^{-1}(\mathbf{DM}(\mathbf{z}^k+\mathbf{v}^k)+\mathbf{DNx}^k- \mathbf{DCB(Qx}^k+\mathbf{u}^k)-\mathbf{x}^k)\rangle\\ =&\langle\mathbf{M}^{*}\mathbf{x}-\mathbf{M}^{*}\mathbf{x
			}^k,\mathbf{z}-\mathbf{(z}^k+\mathbf{v}^k)\rangle+\langle\mathbf{x}-\mathbf{x}^k,(\mathbf{N}-\mathbf{D}^{-1})\mathbf{x}-(\mathbf{N}-\mathbf{D}^{-1})\mathbf{x}^k\rangle\\
			&+\left \langle \mathbf{x}-\mathbf{x}^k,\mathbf{CB(Q}\mathbf{x}^k+\mathbf{u}^k)- \mathbf{CB(Qx)} \right \rangle .
			\endaligned
		\end{equation}
		The first term on the RHS of  (\ref{chang}) can be expressed as
		\begin{equation}\label{101}
			\aligned
			&\langle\mathbf{M}^{*}\mathbf{x}-\mathbf{M}^{*}\mathbf{x}^k
			,\mathbf{z}-\mathbf{(z}^k+\mathbf{v}^k)\rangle \\
			=&\langle \frac1\gamma(\Id-T)(\mathbf{z})-\frac{1}{\gamma_k}(\mathbf{z}^k-\mathbf{z}^{k+1}),
			\mathbf{z}-\mathbf{(z}^k+\mathbf{v}^k)\rangle \\
			=&\frac{1}{\gamma_k}\langle \mathbf{z}^{k+1}-\mathbf{z}^k,
			\mathbf{z}-\mathbf{(z}^k+\mathbf{v}^k)\rangle \\
			=&\frac{1}{2\gamma_k}\left(\|\mathbf{z}^{k+1}-\mathbf{z}^k-
			\mathbf{v}^k\|^{2}+\|\mathbf{z}^k-\mathbf{z}\|^{2}-\|\mathbf{z}^{k+1}-
			\mathbf{z}\|^{2}-\|\mathbf{v}^k\|^{2}\right),
			\endaligned
		\end{equation}
		where the third equality comes from  \eqref{ceq}. The second term on the RHS of  (\ref{chang}) can be rewritten as
		\begin{equation}\label{102}
			\aligned
			&\langle\mathbf{x}  -\mathbf{x}^k,(\mathbf{N}-\mathbf{{D}}^{-1})\mathbf{x}-(\mathbf{N}-\mathbf{{D}}^{-1})\mathbf{x}^k\rangle \\
			=&\frac{1}{2}\langle\mathbf{x}-\mathbf{x}^k,(\mathbf{M}\mathbf{M}^*+2\mathbf{N}-2\mathbf{D}^{-1})
			(\mathbf{x}-{\mathbf{x}}^k)\rangle-\frac{1 }{2}\|\mathbf{M}^*\mathbf{x}-\mathbf{M}^*\mathbf{x}^k\|^2 \\	=&\frac{1}{2}\langle\mathbf{x}-\mathbf{x}^k,(\mathbf{M}\mathbf{M}^{*}+\mathbf{N}+\mathbf{N}^*-2\mathbf{D}^{-1})(\mathbf{x}-\mathbf{x}^k)\rangle -\frac{1}{2\gamma_k^2}\|\mathbf{z}^{k+1}-\mathbf{z}^k\|^2.
			\endaligned
		\end{equation}
		To estimate the last term on the RHS of  (\ref{chang}),
		Young's inequality and \eqref{BB} give
		\begin{equation}\label{103}
			\aligned
			&\langle \mathbf{x}-\mathbf{x}^k,\mathbf{CB(Qx}^k+\mathbf{u}^k)- \mathbf{CB(Qx)}  \rangle \\
			=&- \langle \mathbf{C}^*(\mathbf{x}-\mathbf{x}^k)- [\mathbf{Qx}-\mathbf{(Qx}^k+\mathbf{u}^k)  ],(\mathbf{L}^{-1})^{\frac{1}{2}}\mathbf{L}^{\frac{1}{2}}[\mathbf{B(Qx)-B(Qx}^k+\mathbf{u}^k) ]\rangle \\
			&-\langle \mathbf{Qx}-\mathbf{(Qx}^k+\mathbf{u}^k) ,\mathbf{B(Qx)-B(Qx}^k+\mathbf{u}^k) \rangle\\
			\le&\frac{1}{4} \| (\mathbf{C^*-Q})\mathbf{x}-(\mathbf{C^*-Q})\mathbf{x}^k+\mathbf{u}^k \|_{\mathbf{L}^{-1}} ^2+ \| \mathbf{B(Qx)-B(Qx}^k+\mathbf{u}^k)  \|_{\mathbf{L}} ^2\\
			& - \| \mathbf{B(Qx)-B(Qx}^k+\mathbf{u}^k) \|_{\mathbf{L}} ^2\\
			=&\frac{1}{4}\left \| (\mathbf{C^*-Q})\mathbf{x}-(\mathbf{C^*-Q})\mathbf{x}^k+\mathbf{u}^k\right \|_{\mathbf{L}^{-1}} ^2\\
			\leq&\frac{1}{4}(1+\frac{1}{\theta})\left \| (\mathbf{C^*-Q})(\mathbf{x}-\mathbf{x}^k)\right \|_{\mathbf{L}^{-1}} ^2+\frac{1}{4}(1+\theta)\Vert \mathbf{u}^k \Vert_{\mathbf{L}^{-1}}^2\\
			=& \frac{1}{4}(1+\frac{1}{\theta})\left \langle \mathbf{x-x}^k, (\mathbf{C^*-Q})^*{\mathbf{L}^{-1}}(\mathbf{C^*-Q})(\mathbf{x}-\mathbf{x}^k) \right \rangle+\frac{1}{4}(1+\theta)\Vert \mathbf{u}^k \Vert_{\mathbf{L}^{-1}}^2.
			\endaligned
		\end{equation}
		Substituting (\ref{101}), (\ref{102}) and (\ref{103}) into (\ref{chang}), followed by multiplying by $2\gamma_k$, gives
		\begin{equation*}
			\aligned
			&\| \mathbf{z}^{k+1}-\mathbf{z}  \|^2\\
			\le& \| \mathbf{z}^k-\mathbf{z} \|^2-\frac{1}{\gamma_k }  \|
			\mathbf{z}^{k+1}-\mathbf{z}^k  \|^2 + \| \mathbf{z}^{k+1}-\mathbf{z}^k-\mathbf{v}^k  \|^2- \| \mathbf{v}^k  \|^2+   \frac{\gamma_k(1+\theta)}{2}\Vert \mathbf{u}^k \Vert_{\mathbf{L}^{-1}}^2\\ &+\gamma_k\langle\mathbf{x}-\mathbf{x}^k,[\mathbf{M}\mathbf{M}^{*}+
			\mathbf{N}+\mathbf{N}^*-2\mathbf{D}^{-1}+\frac{1+\theta}{2\theta}
			(\mathbf{C^*-Q})^*{\mathbf{L}^{-1}}(\mathbf{C^*-Q})](\mathbf{x}-\mathbf{x}^k)\rangle \\			=&\|\mathbf{z}^k-\mathbf{z}\|^2+\frac{\gamma_k-1}{\gamma_k}\| \mathbf{z}^{k+1}-\mathbf{z}^k +\frac{\gamma_k}{1-\gamma_k}\mathbf{v}^k \|^2 + \frac{\gamma_k}{1-\gamma_k}\|\mathbf{v}^k\|^2+ \frac{\gamma_k  (1+\theta)}{2}\Vert \mathbf{u}^k \Vert_{\mathbf{L}^{-1}}^2 \\ &+\gamma_k\langle\mathbf{x}-\mathbf{x}^k,[\mathbf{M}\mathbf{M}^{*}+\mathbf{N}+
			\mathbf{N}^*-2\mathbf{D}^{-1}+\frac{1+\theta}{2\theta}
			(\mathbf{C^*-Q})^*{\mathbf{L}^{-1}}(\mathbf{C^*-Q})](\mathbf{x}-\mathbf{x}^k)\rangle.
			\endaligned
		\end{equation*}
		According to (b) of Remark \ref{rem11}(2) and Assumption \ref{assumption2}(ii), we obtain
		\begin{equation*}
			\aligned
			\| \mathbf{z}^{k+1}-\mathbf{z} \|^2&+\frac{1-\gamma_k}{\gamma_k}\| \mathbf{z}^{k+1}-\mathbf{z}^k +\frac{\gamma_k}{1-\gamma_k}\mathbf{v}^k \|^2\\
			&\le  \| \mathbf{z}^k-\mathbf{z}  \|^2+\frac{\gamma_k}{1-\gamma_k}\|\mathbf{v}^k\|^2+\frac{\gamma_k (1+\theta )}{2} \| \mathbf{u}^k  \|_{\mathbf{L}^{-1}} ^2.
			\endaligned
		\end{equation*}
		From (\ref{ln}), we know the inequality (\ref{lem1}) holds. Finally, (\ref{dfd}) follows from (\ref{tiaojian}).
	\end{proof}
	
	The following theorem is our main result regarding the convergence of  Algorithm  \ref{Alg12}.
	\begin{theorem}\label{theo}
		{
			Let  $F_1,\dots,F_n: \mathcal{H}\rightrightarrows \mathcal{H}$ be maximally monotone, and let  $B_i: \mathcal{H}\to \mathcal{H} $ be $\frac{1}{L_i}$-cocoercive  for $i=1,\dots,m$ with  $\zer({\textstyle \sum_{i=1}^{n}}  F_i+{\textstyle \sum_{i=1}^{m}} B_i)\neq \emptyset $. Suppose Assumptions  \ref{assumption2} and \ref{asm11}   hold.	 Let $\{\mathbf{x}^k\}\subseteq\mathcal{H}^{n}$ and $\{\mathbf{z}^k\}\subseteq\mathcal{H}^{n-1}$  be generated by  Algorithm \ref{Alg12}, and $\{\mathbf{u}^k\}\subseteq\mathcal{H}^{m}$ and $\{\mathbf{v}^k\}\subseteq\mathcal{H}^{n-1}$  satisfy \eqref{tiaojian}. Then the following assertions hold:}
		\begin{itemize}
			\item[{\rm(i)}]  The sequence $\{l^2_k\}$ is summable and the sequences   $\{\mathbf{u}^k\}$ and $\{\mathbf{v}^k\}$ converge to zero.
			\item[{\rm(ii)}] For all  $\mathbf{z}\in \Fix T$, the sequence $\{\|\mathbf{z}^{k}-\mathbf{z}  \| \}$ converges.
			\item[{\rm(iii)}]  As $k\to \infty $, we have $\mathbf{z}^{k+1}-\mathbf{z}^k\to 0$ and $ {\textstyle \sum_{i=1}^{n}}s_ix_i^k\to 0 $ for all $s\in \mathbb{R}^n$ with $ {\textstyle \sum_{i=1}^{n}}s_i= 0 $.
			\item[{\rm(iv)}]  The sequence $\{\mathbf{z}^k\}$ converges weakly to a point in $\Fix T$.
			\item[{\rm(v)}] The sequence $\{\mathbf{x}^k\}$ converges weakly to a point $(\bar{x},\ldots,\bar{x})\in\mathcal{H}^n$ such that \vskip 1mm $\bar{x}\in\zer\left ( \sum_{i=1}^{n}F_i+\sum_{i=1}^{m}B_i   \right )$.
			
		\end{itemize}					
	\end{theorem}
	\begin{proof}
		{\rm(i)}
		Now, we start by proving (i) via (\ref{dfd}). Let $\mathbf{z}\in \Fix T$ and  $N\in\mathbb{N}_+$. We sum (\ref{dfd}) for $k=1,2,\ldots,N$ to obtain
		$$ \| \mathbf{z}^{N+1}-\mathbf{z}  \|^2+l_N^2+\sum_{k=1}^{N-1}(1-\xi _k)l^2_k \le  \| \mathbf{z}^{1}-\mathbf{z}  \|^2+\xi _0l_0^2. $$
		Then, rearranging the terms gives
		\begin{equation}\label{24}
			\aligned
			\sum_{k=1}^{N} (1-\xi _k)l^2_k&\le  \| \mathbf{z}^1-\mathbf{z} \|^2- \| \mathbf{z}^{N+1}-\mathbf{z}  \|^2-\xi _Nl^2_N+\xi_0l^2_0\\
			& \le  \| \mathbf{z}^1-\mathbf{z}  \|^2+\xi_0l^2_0 .
			\endaligned
		\end{equation}
		Since the RHS of \eqref{24} is independent of $N$, we conclude that
		$\sum_{k=1}^{\infty }(1-\xi _k)l^2_k< \infty . $
		Thanks to $\{\xi _k\}\subseteq[0, 1-\epsilon]$ from Assumption \ref{assumption2}(iii), we have
		$\sum_{k=1}^{\infty }l^2_k< \infty , $
		which implies
		\begin{equation}
			\label{ln2}
			\lim_{k \to \infty}l^2_k=0.
		\end{equation}
		Then, from (\ref{tiaojian}), Assumption \ref{assumption2}(i) and (ii) we deduce that $\lim_{k \to \infty}  \mathbf{u}^k= 0$ and $\lim_{k \to \infty} \mathbf{v}^k=0$.
		
		(ii)   From (\ref{dfd}) and Assumption \ref{assumption2}(iii), we deduce
		that the sequence $\{ \| \mathbf{z}^k-\mathbf{z}  \|^2+l_{k-1}^2 \}$ is nonincreasing and hence convergent. Furthermore, from (\ref{ln2}),  the sequence $\{ \| \mathbf{z}^k-\mathbf{z}  \|\}$ converges, which implies that the sequence $\{\mathbf{z}^k\}$ is bounded.
		
		(iii)	 From (\ref{ln}), (\ref{ln2}), Assumption \ref{assumption2}(ii) and $\lim_{k \to \infty} \mathbf{v}^k=0$, we deduce that $\lim_{k \to \infty}  ( \mathbf{z}^{k+1}-\mathbf{z}^k )=0.$  According to Assumption \ref{asm11}(a),  let $s\in\{a\in\mathbb{R}^n:\sum_{i=1}^na_i=0\}=(\ker M^{\top})^\perp=\operatorname{range}M.$ Then, there exists $v\in \mathbb{R}^{n-1}$ such that $s=-Mv$ and hence
		\begin{equation*}
			\aligned
			\sum_{i=1}^ns_ix_i^k & =(s^{\top}\otimes\Id)\mathbf{x}^{k}=-\left((v^{\top}M^{\top})\otimes\Id\right)\mathbf{x}^{k} \\
			&=-(v^{\top}\otimes\Id)\mathbf{M^*x}^{k}=\frac{1}{\gamma_k}(v^{\top}\otimes\Id)(\mathbf{z}^{k+1}-\mathbf{z}^{k})\to0.
			\endaligned
		\end{equation*}

		(iv) and (v)  Let $\mathbf{x}^k=J_{ \mathbf{DF}}(\mathbf{y}^k)$ where $\mathbf{y}^k=\mathbf{DM}(\mathbf{z}^k+\mathbf{v}^k)+\mathbf{DN}\mathbf{x}^k- \mathbf{DCB(Q}\mathbf{x}^k+\mathbf{u}^k)$. We claim that the sequence $\{\mathbf{x}^k\}$ is bounded. To see this, let ${d_1}M_1$ denote the first row of the matrix $DM$. By Assumption \ref{asm11}(b) and Remark \ref{rem11}(2), we have
		$$
		x_1^k=J_{{d_1}F_1}(y_1^k)=J_{{d_1}F_1}({d_1}
		\mathbf{M}_1(\mathbf{z}^k+\mathbf{v}^k)).$$
		From the nonexpansivity of resolvents  and boundedness of $\{\mathbf{z}^k\}$ and $\{\mathbf{v}^k\}$, it  follows that $\{x_1^k\}$ is also bounded. By (iii),  $\{\mathbf{x}^k\}$ is bounded, as claimed. Let $ \mathbf{\bar{z}}=(\bar{z}_1,\dots,\bar{z}_{n-1})\in\mathcal{H}^{n-1}$ be an arbitrary weak cluster point of $\{\mathbf{z}^k\}$. Then there exists a point $\mathbf{\bar{x}}\in\mathcal{H}^n$ such that
		$(\mathbf{\bar{z}},\mathbf{\bar{x}})$ is a weak cluster point of $\{(\mathbf{z}^k,\mathbf{x}^k)\}$, where $\mathbf{\bar{x}}=(\bar{x},\dots,\bar{x})$ according to  (iii).
		Denote $\mathbf{w}^k=\mathbf{DM}(\mathbf{z}^k+\mathbf{v}^k)+\mathbf{DN}\mathbf{x}^k$. Then $\mathbf{y}^k=\mathbf{w}^k- \mathbf{DCB(Q}\mathbf{x}^k+\mathbf{u}^k)$ and $\mathbf{w}=\mathbf{DM}\mathbf{\bar{z}}+\mathbf{DN\bar{x}}$ is a weak cluster point of $\{\mathbf{w}^k\}$. Using  $\range {M}=\{a\in\mathbb{R}^n:\sum_{i=1}^na_i=0\}$, we deduce ${\textstyle \sum_{i=1}^{n}}\frac{w_i^k}{d_i}= {\textstyle \sum_{i,j=1}^{n}} N_{ij}x_j^k$ from $\mathbf{D}^{-1}\mathbf{w}^k=\mathbf{M}(\mathbf{z}^k+\mathbf{v}^k)+\mathbf{N}\mathbf{x}^k$. Then, we get
		\begin{equation}\label{FT}
			\aligned
			\Phi\begin{pmatrix}\frac{1}{d_1}( w_1^k-x_1^k)\\\frac{1}{d_2}(w_2^k-x_2^k)+ b_2^k \\ \vdots \\\frac{1}{d_{n-1}}( w_{n-1}^k-x_{n-1}^k)+ b_{n-1}^k\\x_n^k\end{pmatrix} &\ni \begin{pmatrix}x_1^k-x_n^k \\ x_2^k-x_n^k\\ \vdots \\ x_{n-1}^k-x_n^k\\ {\textstyle \sum_{i=1}^{n}}\frac{1}{d_i} (w_i^k-x_i^k) + {\textstyle \sum_{i=2}^{n}}b_i^k \end{pmatrix}\\
			&=\begin{pmatrix} x_1^k-x_n^k \\ x_2^k-x_n^k\\ \vdots \\ x_{n-1}^k-x_n^k\\ {\textstyle \sum_{i,j=1}^{n}}N_{ij}x_j^k-{\textstyle \sum_{i=1}^{n}}\frac{1}{d_i} x_i^k +  {\textstyle \sum_{i=2}^{n}}b_i^k\end{pmatrix},
			\endaligned
		\end{equation}
		where $b_i^k={\textstyle \sum_{j=1}^{m}}C_{ij} B_{j}(x_i^k)-{\textstyle \sum_{j=1}^{m}}C_{ij} B_{j}( {\textstyle \sum_{h=1}^{i-1}}Q_{jh}x_h^k+u_{j}^k) $, $i=2,\ldots,n$ and the operator $\Phi: \mathcal{H}^n\rightrightarrows \mathcal{H}^n$ is defined by
		\begin{equation}\label{phi}
			\Phi= \begin{pmatrix}
				F_1^{-1}\\
				(F_2+{\textstyle \sum_{j=1}^{m}}C_{2j} B_j)^{-1}\\
				\vdots \\
				(F_{n-1}+{\textstyle \sum_{j=1}^{m}}C_{n-1,j} B_{j})^{-1}\\
				F_n+{\textstyle \sum_{j=1}^{m}}C_{nj} B_{j}
			\end{pmatrix}+\begin{pmatrix}
				0 & 0 & 0  &\cdots  &0  & -\Id\\
				0 & 0 & 0 &\cdots  &0  & -\Id \\
				\vdots&\vdots   & \ddots  &\vdots   &\vdots  \\
				0 & 0 & 0  &\cdots  &0  & -\Id \\
				\Id& \Id &\Id&\cdots & \Id &0
			\end{pmatrix}.
		\end{equation}
		According to  (i), (iii), Assumption \ref{asm11}(b)  and the Lipschitz continuity  of $B_i$ we have $\lim_{k \to \infty} b_i^k=0$ for $i=2,\ldots,n.$
		As the sum of two maximally monotone operators is again maximally monotone provided that one of the operators has full domain \cite[Corollary 24.4 (i)]{BC2011}, it follows that $\Phi$ is maximally monotone. Consequently, its graph is sequentially closed in the weak-strong topology  \cite[Proposition 20.32]{BC2011}.  Note also that the RHS of (\ref{FT}) converges strongly to zero as a consequence of (i), (iii) and (a) of Remark \ref{rem11}(2).	
		Taking the limit  along a subsequence of $\{(\mathbf{z}^k,\mathbf{x}^k)\}$ which converges weakly to $(\mathbf{\bar{z}},\mathbf{\bar{x}})$ in (\ref{FT}), and using the weak-strong topology of $\Phi$, we obtain
		$$\Phi
		\begin{pmatrix}
			\frac{1}{d_1}(w_1-\bar{x})\\
			\vdots \\
			\frac{1}{d_{n-1}}(w_{n-1}-\bar{x})\\
			\bar{x}
		\end{pmatrix}\ni \begin{pmatrix}
			0 \\
			\vdots \\
			0 \\
			0
		\end{pmatrix}.$$ Then, it follows
		\begin{equation}\label{phi0}
			\left\{\begin{array}{l}
				F_{1}(\bar{x}) \ni \frac{1}{d_1}(w_1-\bar{x} ),\\
				(F_{i}+{\textstyle \sum_{j=1}^{m}}C_{ij}B_{j})(\bar{x})\ni \frac{1}{d_i}(w_i-\bar{x}),\ \forall i \in[2, n-1],\\
				(F_{n}+{\textstyle \sum_{j=1}^{m}}C_{nj}B_{j})(\bar{x}) \ni-\sum_{i=1}^{n-1}\frac{1}{d_i}(w_{i}-\bar{x}),
			\end{array}\right.
		\end{equation}
		which with Lemma \ref{lem11}(c) yields $\mathbf{\bar{z}}\in\Fix T$. Hence $\omega_w(\mathbf{{z}}^k)\subseteq\Fix(T)$. Using (ii) and  Lemma \ref{lle}, it follows that $\{\mathbf{z}^k\}$ converges  weakly to a point in $\Fix T$.  	From \eqref{phi0}, we also have	  $\bar{x}= J_{{d_1}F_1}(w_1)=J_{{d_1}F_1}({d_1}\mathbf{M}_1\mathbf{\bar{z}})$ and $\bar{x}\in\zer\left ( \sum_{i=1}^{n}F_i+\sum_{i=1}^{m}B_i   \right )$. Therefore $\mathbf{\bar{x}}=(\bar{x},\dots,\bar{x})$  is the unique weak sequential cluster point of the  sequence $\{\mathbf{x}^k\}$ and  $\{\mathbf{x}^k\}$ converges  weakly to $\mathbf{\bar{x}}$.
	\end{proof}

	\section{Numerical experiments}
	Numerical experiments in \cite{Anton} illustrate that among some schemes in  Algorithm \ref{AlgG}, the split-forward-backward+ (SFB+) and the adapted graph forward-backward splitting (aGFB)  achieve favorable numerical results. Therefore, in this section, we consider split-forward-backward+ with deviations (dev\_SFB+) and the adapted graph forward-backward splitting with deviations (dev\_aGFB).

	\subsection{Problem description}
	A classic problem for testing optimization algorithms is the Markowitz portfolio optimization  in \cite{DRS}. Here we focus on this problem of the form:
	\begin{equation}\label{mark}
		\underset{x \in \Delta}{\min} \, \frac{1}{2}x^\top\Lambda  x - r^\top x + \frac{\delta}{2} \| x \|^2 + \sum_{i=1}^n \bigl| x_i - (x_0)_i \bigr| + \sum_{i=1}^n \bigl| x_i - (x_0)_i \bigr|^{3/2},
	\end{equation}
	where $\delta>0,r\in \mathbb{R}^n$ is a vector of estimated assets returns, $\Lambda\in \mathbb{R}^{n\times n} $ is the estimated
	covariance matrix of returns (which is  symmetric positive semidefinite), $\Delta =\{x\in\mathbb{R}^n:x_i\ge 0,\forall i=1,\ldots,n, x_1+\dots +x_n=1\}$ is the standard simplex, and
	$x_0 \in \mathbb{R}^n$ is the initial position.
	We introduce  the indicator function $\mathbb{I}_\Delta$ of the constraint set $\Delta$ to reformulate \eqref{mark} as an unconstrained optimization problem and	then  split the objective function into the following form:
	\begin{equation}\label{wenti}
		\min_{x \in \mathbb{R}^n} \, \underbrace{\frac{1}{2}x^\top \Lambda x- r^\top x}_{f_1(x)} +
		\underbrace{ \frac{\delta}{2} \| x \|^2 }_{f_2(x)}
		+ \underbrace{\sum_{i=1}^n \vert x_i - (x_0)_i \vert}_{g_1(x)}
		+ \underbrace{\sum_{i=1}^n \vert x_i - (x_0)_i \vert^{3/2}}_{g_2(x)}
		+ \underbrace{\mathbb{I}_\Delta(x)}_{g_3(x)}.
	\end{equation}
	%We denote by $\alpha$ the largest eigenvalue of $\Lambda $, i.e., the Lipschitz constant of $w\mapsto \Lambda w$.
	Let $F_i=\partial g_i$ for $i=1,2,3$ and $B_i= \nabla f_i$ for $i=1,2$. Then, the optimization problem \eqref{wenti} is equivalent to  the inclusion problem: find $x \in \mathbb{R}^n$ such that
	\begin{equation*}
		0 \in (F_1+F_2+F_3+B_1+B_2)(x),
	\end{equation*}
	where $F_i$, $i=1,2,3$ are the maximally monotone operators, as shown in  \cite[Theorem 20.48]{BC2011}. By \cite[Corollary 18.17]{BC2011}, $B_i$, $i=1,2$ are $\frac{1}{L_i}$-cocoercive, where $L_1$ is the largest eigenvalue of $\Lambda$ and $L_2=\delta$.

	\subsection{Experiment settings and results}
	In the experiment, the estimated returns $r$ and the matrix $\Lambda$ are derived from real data. Details are available in $\hbox{data}$\footnote{Code and data can be downloaded from https://github.com/TraDE-OPT/portopt-siopt21.} for  $200$ market days, and $n = 53$ assets. We set $\delta=6$. Consider the following two cases of the problem \eqref{mark}:
	
	\noindent \textbf{Case 1}: the initial state $x_0$ is chosen randomly; \\
	\noindent  \textbf{Case 2}: $x_0$ is the output of Case 1 for the same optimization problem but $20$ days later. 
	
	We test aGFB and SFB+ and their formulas with deviations to show the effectiveness of the deviations. For the four algorithms, we take the initial point $z^0_i=(0,\dots,0)\in \mathbb{R}^n$, $i=1,2$, and  choose the matrices $M$, $C$ and $Q$  as in \cite{Anton}.  We take $\gamma_k=0.9$ for all algorithms.  For dev\_SFB+ and dev\_aGFB,  $\theta$  is tuned manually to a nearly optimal value. Initialize $u_i^0=(0,\dots,0)\in \mathbb{R}^n$ and  $v_i^0=(0,\dots,0)\in \mathbb{R}^n$  for $i=1,2$. The algorithms terminate  when
	$$
	\|x_3^k-x^*\|<10^{-8}
	$$
	where $x^*$ is the exact solution.
	
	There are some ways for selecting deviation vectors that satisfy the deviation condition. In order to select relatively optimal deviation vectors, we refer to \cite{Ban} and adopt a three-layer convolutional neural network model. 
	
\begin{figure}[H]
		\centering  % 整个图居中
		% 第1行：2个子图
		\begin{subfigure}[b]{0.42\columnwidth}  
			\includegraphics[width=\textwidth]{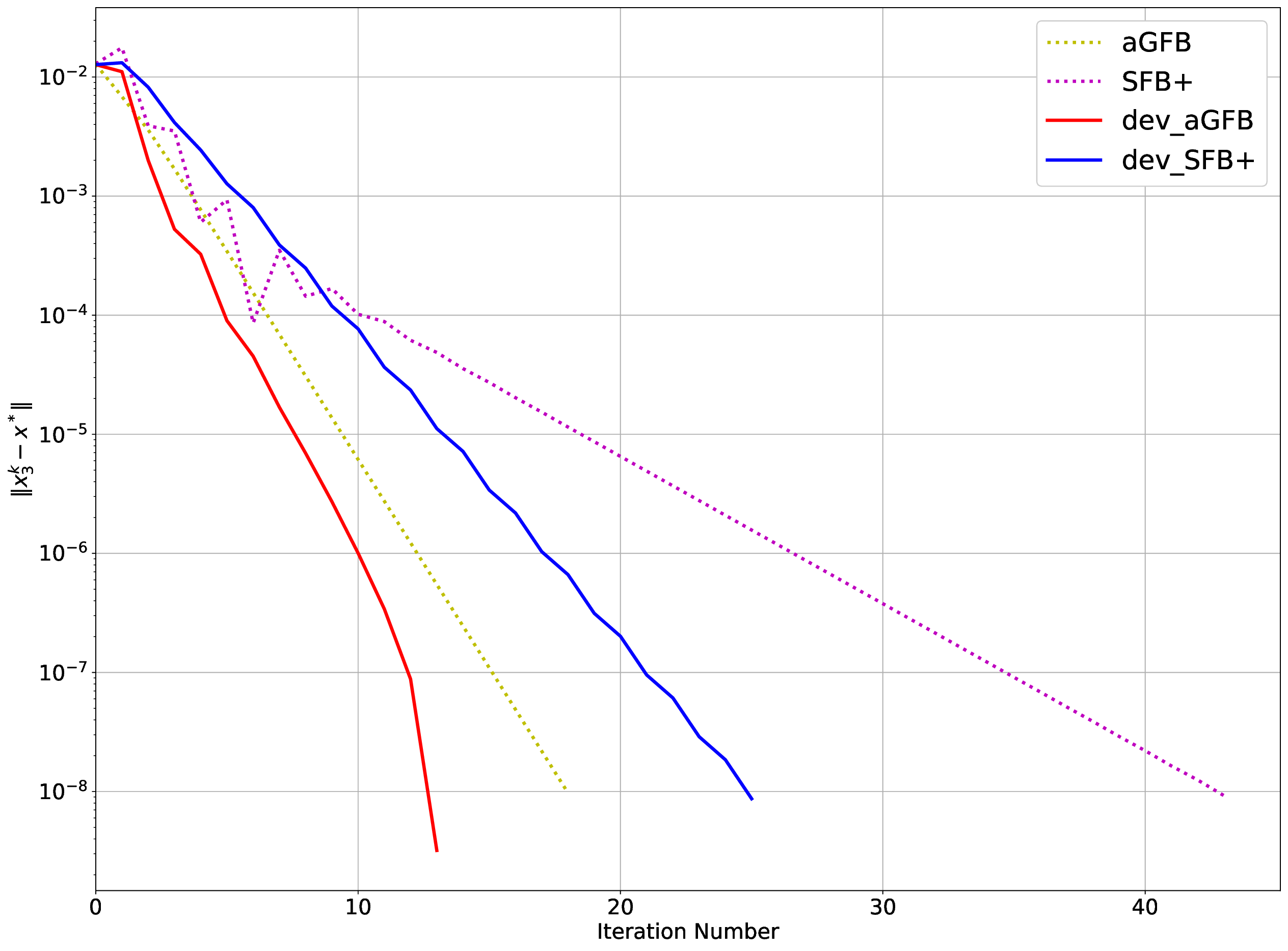}  
			\caption{Case 1}
		\end{subfigure}
		\begin{subfigure}[b]{0.42\columnwidth}
			\includegraphics[width=\textwidth]{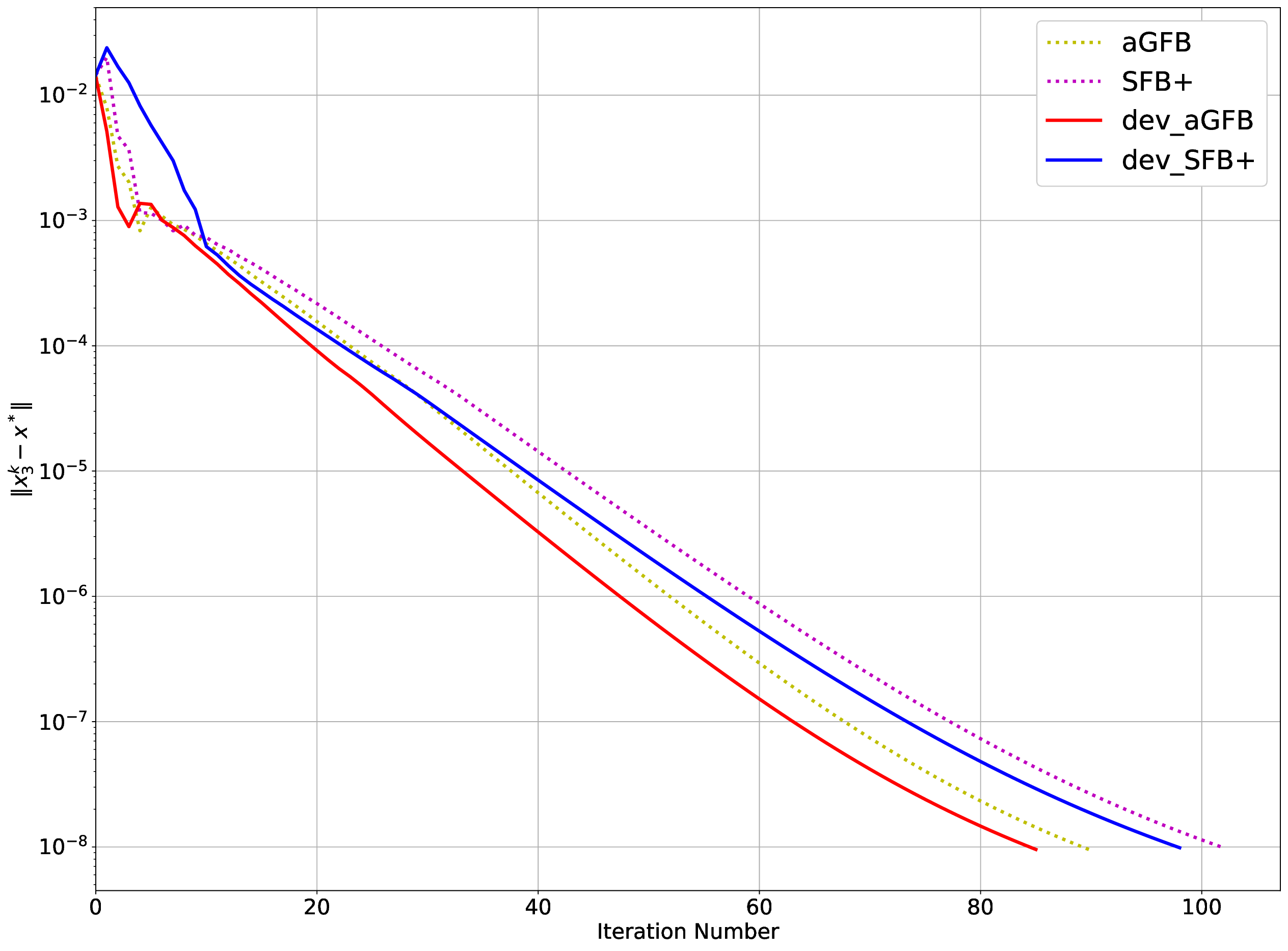}
			\caption{Case 2}
		\end{subfigure}
		% 主图标题
		\caption{The decay of the errors  of four algorithms with the numbers of the iterations  under Case 1 and Case 2.}
		\label{fig}
	\end{figure}

	Figure 	\ref{fig} presents the decay of the errors  of four algorithms with the numbers of the iterations. It demonstrates the improvement  in the iterative methods achieved by the appropriate selection of the deviations.
	We  randomly choose $x_0$ and  report  arithmetical averages of the required
	numbers of iterations taken to run the algorithms 50 times in Table \ref{table}. 
	It can be observed from Figure 	\ref{fig} and Table \ref{table} that for Case 1, the dev\_aGFB and dev\_SFB+ algorithms perform significantly better than their counterparts (aGFB and SFB+) while for Case 2, the improvements of dev\_aGFB and dev\_SFB+ on aGFB and SFB+ are slightly. Furthermore, dev\_aGFB behaves best among four algorithms. 
	
	\begin{table}[htbp]
		\centering
		\caption{Comparisons of the numbers of iterations of four algorithms.}
		\label{tab:numerical-results-transposed}
		\begin{tabular}{lccccc} 
			\toprule
			& \multicolumn{4}{c}{Algorithms} \\
			\cmidrule(lr){2-5}
			& aGFB & SFB+ & dev\_aGFB & dev\_SFB+  \\ 
			\midrule
			Case 1  & 28 & 44 & 23.64 & 26.34  \\ 
			Case 2  & 89.06 & 91.16 & 85.12 & 87.56  \\
			\bottomrule
		\end{tabular}
		\label{table}
	\end{table}

\section*{Acknowledgements}
The third author was supported in part by National Natural Science Foundation of China (No. 12271273).


\vskip 2mm	
	
	\begin{thebibliography}{00}
		
		\bibitem{ring}
		Arag\'{o}n-Artacho, F.J., Malitsky, Y., Tam, M.K. et al. Distributed forward-backward methods for ring networks. Computational Optimization and Applications. 86, 845-870 (2023).
		
		
		\bibitem{Moursi}
		Moursi, W.M. The forward-backward algorithm and the normal problem.  Journal of Optimization Theory and Applications. 176, 605-624 (2018).
		
		\bibitem{Chen}
		Chen, G.H., Rockafellar, R.T. Convergence rates in forward-backward splitting. SIAM Journal on Optimization. 7(2), 421-444 (1997).	
		
		\bibitem{Ryu}
		Ryu, E.K.  Uniqueness of DRS as the 2 operator resolvent-splitting and impossibility of 3 operator resolvent-splitting. Mathematical Programming. 182, 233-273 (2020).
		
		\bibitem{fran}
		Arag\'{o}n-Artacho, F.J.,  Campoy, R., L\'{o}pez-Pastor, C. Forward-backward algorithms devised by graphs. https://arxiv.org/abs/2406.03309.
		
		\bibitem{Frugal}
		Morin, M., Banert, S., Giselsson, P. Frugal splitting operators: representation, minimal lifting, and convergence. SIAM Journal on Optimization. 34(2), 1187-2168 (2024).
		
		\bibitem{Bred}
		Bredies, K., Chenchene, E., Lorenz, D.A. et al. Degenerate preconditioned proximal point algorithms. SIAM Journal on Optimization. 32(3), 1499-2459 (2022).	
		
		\bibitem{Malit}
		Malitsky, Y., Tam, M.K. Resolvent splitting for sums of monotone operators with minimal lifting. Mathematical Programming.  201(1), 231-262 (2023).
		
		\bibitem{Pierra}
		Pierra, G. Decomposition through formalization in a product space. Mathematical Programming. 28, 96-115 (1984).
		
		\bibitem{Bre}
		Bredies, K., Chenchene, E., Naldi, E. Graph and distributed extensions of the Douglas--Rachford method. SIAM Journal on Optimization. 34(2), 1569-1594 (2024).
		
		\bibitem{Anton}
		\AA kerman, A.,  Chenchene, E., Giselsson, P. et al. Splitting the forward-backward algorithm: a full characterization. 	https://arxiv.org/abs/2504.10999.
		
		\bibitem{minh}
		Dao, M.N., Tam, M.K., Truong, T.D. A general approach to distributed operator splitting. https://arxiv.org/abs/2504.14987.
		
		
		\bibitem{Sedeghi}
		Sadeghi, H., Banert, S., Giselsson, P. Forward-backward splitting with deviations for monotone inclusions. Applied Set-Valued Analysis and Optimization. 6(2), 113-135 (2024).
		
		\bibitem{L2O}
		Chen, T., Chen, X., Chen, W. et al. Learning to optimize: a primer and a benchmark. The Journal of Machine Learning Research. 23(1), 8562-8620 (2022).
		
		
		\bibitem{Ban}
		Banert, S., Rudzusika, J., \"Oktem, O. et al. Accelerated forward-backward optimization using deep learning. SIAM Journal on Optimization. 34(2), 1236-1263 (2024).
		
		\bibitem{Qin}
		Qin, L., Huang, X., Dong, Q.L. et al. Accelerated Douglas--Rachford splitting algorithm using neural network. preprint.
		
		\bibitem{Hu-Davis-Yin}
		Hu, Z., Dong, Q.L. A three-operator splitting algorithm with deviations for generalized DC programming. Applied Numerical Mathematics. 191, 62-74 (2023).
		
		\bibitem{Qin-FBHF}
		Qin, L., Dong, Q.L., Zhang, Y. et al. Forward-backward-half forward splitting algorithm with deviations. Optimization. 74(9), 2137-2158 (2024).
		
		\bibitem{BC2011}
		Bauschke, H.H., Combettes, P.L. \textit{Convex Analysis and Monotone Operator Theory in Hilbert Spaces}, 2nd ed. Springer, New York, (2017).
		
		\bibitem{Matt}
		Tam, M.K. Frugal and decentraised resolvent splitting defined by nonexpansive operators. Optimization Letters. 18(7), 1-19 (2023).	
		
		\bibitem{DRS}
		Brodie, J., Daubechies, I., De Mol, C. et al. Sparse and stable Markowitz portfolios. Proceedings of the National Academy of Sciences of the United States of America. 106(30), 12267-12272 (2009).
		
		
		
		
		
		
		
		
		
		
		
		
		
		
		
		
		
		
		
		%\bibitem{Bruck}
		%Ronald E., Bruck Jr. An iterative solution of a variational inequality for certain monotone operators in Hilbert space. Bulletin of the American Mathematical Society. 81, 890-892 (1975).
		
		%\bibitem{Lions}
		%Lions, P.L., Mercier, B. Splitting algorithms for the sum of two nonlinear operators. SIAM Journal on Numerical Analysis.  16(6), 964-979 (1979).
		
		%\bibitem{Passty}
		%Passty, G.B. Ergodic convergence to a zero of the sum of monotone operators in Hilbert space. Journal of Mathematical Analysis and Applications. 72(2), 383-390 (1979).
		
		
		
		%\bibitem{Karol}
		%Karol, G., Yann, L. Learning fast approximations of sparse coding. Proceedings of the 27th International Conference on International Conference on Machine Learning. 399-406 (2010).
		
		%\bibitem{Banert}
		%Banert, S., Ringh, A., Adler, J. et al. Data-driven nonsmooth optimization. SIAM Journal on Optimization. 30(1), 102-131 (2020).
		
		
		%\bibitem{Douglas}
		%Douglas, J., Rachford, H. On the numerical solution of heat conduction problems in two and three space variables. Transactions of the American Mathematical Society. 82(2), 421-439 (1956).	
		
		%\bibitem{DY}
		%Davis, D., Yin, W. A three-operator splitting scheme and its optimization applications. Set-valued and variational analysis. 25(4), 829-858 (2017).
		
		%\bibitem{Condat}
		%Condat, L., Kitahara, D., Contreras, A. et al. Proximal splitting algorithms for convex optimization: a tour of recent advances, with new twists. SIAM Review. 65(2), 375-435 (2023).
		
		%\bibitem{Rock}
		%Rockafellar, R.T. Monotone operators and the proximal point algorithm. Siam Journal on Control and Optimization.  14(5) 877-898 (1976).
		
		%\bibitem{HeB}
		%He, B.S., Yuan, X. Convergence analysis of primal-dual algorithms for total variation image restoration. SIAM Journal on Imaging Sciences. 5(1), 119-149 (2010).
		
		
		%\bibitem{FBHF}
		%Brice\~{n}o-Arias, L.M., Davis, D. Forward-backward-half forward algorithm for solving monotone inclusions. SIAM Journal on Optimization.  28(4), 2839-2871 (2018).
		
		
		
		
	\end{thebibliography}
\end{document}